\documentclass{notices}
\usepackage{amsmath,amssymb,amsthm}
\usepackage{bm}
\usepackage{graphicx}
\usepackage{hyperref}
\usepackage{xcolor}
\usepackage[shortlabels]{enumitem}
\usepackage{tikz-cd}

\DeclareMathOperator{\ran}{ran}
\DeclareMathOperator{\spn}{span}
\DeclareMathOperator{\Id}{Id}
\DeclareMathOperator{\Real}{Re}
\xdefinecolor{dgred}{rgb}{1,0.5020,0}
\xdefinecolor{dgblue}{rgb}{.2, .2, .702}

\theoremstyle{definition}
\newtheorem{prob}{Problem}
\newtheorem{exmpl}{Example}

\allowdisplaybreaks

\title{Bridging data science and dynamical systems theory}

\author{
  Tyrus Berry
  \affil{
      Tyrus Berry is an assistant professor at George Mason University. His email address is tberry@gmu.edu.
    }
  \and
  Dimitrios Giannakis
  \affil{
      Dimitrios Giannakis is an associate professor at New York University. His email address is dimitris@cims.nyu.edu.}
  \and
  John Harlim
  \affil{ 
      John Harlim is a professor at the Pennsylvania State University. His email address is jharlim@psu.edu.}
}

\begin{document}

\maketitle

Modern science is undergoing what might arguably be called a ``data revolution'', manifested by a rapid growth of observed and simulated data from complex systems, as well as vigorous research on in mathematical and computational frameworks for data analysis. In many scientific branches, these efforts have led to the creation of statistical models of complex systems that match or exceed the skill of first-principles models. Yet, despite these successes, statistical models are oftentimes treated as black boxes, providing limited guarantees about stability and convergence as the amount of training data increases. Black-box models also offer limited insights about the operating mechanisms (physics), the understanding of which is central to the advancement of science. 

In this short review, we describe mathematical techniques for statistical analysis and prediction of time-evolving phenomena, ranging from simple examples such as an oscillator, to highly complex systems such as the turbulent motion of the Earth's atmosphere, the folding of proteins, and the evolution of species populations in an ecosystem. Our main thesis is that combining ideas from the theory of dynamical systems with learning theory provides an effective route to data-driven models of complex systems, with refinable predictions as the amount of training data increases, and physical interpretability through discovery of coherent patterns around which the dynamics is organized. Our article thus serves as an invitation to explore ideas at the interface of the two fields. 

This is a vast subject, and invariably a number of important developments in areas such as deep learning \cite{EEtAl17}, reservoir computing \cites{PathakEtAl18,VlachasEtAl20}, control \cites{KordaMezic18,KlusEtAl20}, and non-autonomous/stochastic systems \cites{Froyland13,KlusEtAl18} are not discussed here.

\section*{Statistical forecasting and coherent pattern extraction}

Consider a dynamical system of the form $ \Phi^t : \Omega \to \Omega$, where $\Omega$ is the state space and $\Phi^t$, $ t \in \mathbb R$, the flow map. For example, $\Omega$ could be Euclidean space $\mathbb R^d$, or a more general manifold, and $\Phi^t$ the solution map for a system of ODEs defined on $\Omega$. Alternatively, in a PDE setting, $\Omega$ could be an infinite-dimensional function space and $ \Phi^t $ an evolution group acting on it. We consider that $\Omega$ has the structure of a metric space equipped with its Borel $\sigma$-algebra, playing the role of an event space, with measurable functions on $\Omega$ acting as random variables, called \emph{observables}. 

In a statistical modeling scenario, we consider that available to us are time series of various such observables, sampled along a dynamical trajectory which we will treat as being unknown. Specifically, we assume that we have access to two observables, $ X : \Omega \to \mathcal X$ and $ Y : \Omega \to \mathcal Y$, respectively referred to as covariate and response functions, together with corresponding time series $ x_0, x_1, \ldots, x_{N-1}$ and $ y_0, y_1, \ldots, y_{N-1}$, where $ x_n = X(\omega_n) $, $ y_n = Y(\omega_n)$, and $ \omega_n = \Phi^{n \, \Delta t}(\omega_0)$. Here, $\mathcal X$ and $\mathcal Y$ are metric spaces, $\Delta t$ is a positive sampling interval, and $\omega_0$ an arbitrary point in $\Omega$ initializing the trajectory. We shall refer to  the collection $ \{ ( x_0, y_0 ), \ldots, ( x_{N-1}, y_{N-1}) \} $ as the training data. We require that $\mathcal Y$ is a Banach space (so that one can talk about expectations and other functionals applied to $Y$), but allow the covariate space $\mathcal X$ to be nonlinear. 

Many problems in statistical modeling of dynamical systems can be expressed in this framework. For instance, in a low-dimensional ODE setting, $X$ and $Y$ could both be the identity map on $\Omega = \mathbb R^d$, and the task could be to build a model for the evolution of the full system state. Weather forecasting is a classical high-dimensional application, where $\Omega$ is the abstract state space of the climate system, and $X$ a (highly non-invertible)  map representing  measurements from satellites, meteorological stations, and other sensors available to a forecaster. The response $Y$ could be temperature at a specific location, $\mathcal Y = \mathbb R$, illustrating that the response space may be of considerably lower dimension than the covariate space. In other cases, e.g., forecasting the temperature field over a geographical region, $\mathcal Y$ may be a function space.  The two primary questions that will concern us here are:
\begin{prob}[Statistical forecasting] \label{probForecast} Given the training data, construct (``learn'') a function $ Z_t : \mathcal X \to \mathcal Y$ that predicts $Y$ at a lead time $ t \geq 0$. That is, $Z_t $ should have the property that $ Z_t \circ X $ is closest to $ Y \circ \Phi^t $ among all functions in a suitable class. 
\end{prob}
\begin{prob}[Coherent pattern extraction] \label{probSpectral} Given the training data, identify a collection of observables $ z_j : \Omega \to \mathcal Y $ which have the property of evolving coherently under the dynamics. By that, we mean that $ z_j \circ \Phi^t$ should be relatable to $z_j$ in a natural way. 
\end{prob}

These problems  have an extensive history of study from an interdisciplinary perspective spanning mathematics, statistics, physics, and many other fields. Here, our focus will be on \emph{nonparametric methods}, which do not employ explicit parametric models for the dynamics. Instead, they use universal structural properties of dynamical systems to inform the design of data analysis techniques. From a learning standpoint, Problems~\ref{probForecast} and~\ref{probSpectral} can be thought of as \emph{supervised} and \emph{unsupervised} learning, respectively. A mathematical requirement we will impose to methods addressing either problem is that they  have a well-defined notion of convergence, i.e., they are refinable, as the number $N$ of training samples increases.

\section*{Analog and POD approaches}

Among the earliest examples of nonparametric forecasting techniques is Lorenz's analog method \cite{Lorenz69b}. This simple, elegant approach makes predictions by tracking the evolution of the response along a dynamical trajectory in the training data (the analogs). Good analogs are selected according to a measure of geometrical similarity between the covariate variable observed at forecast initialization and the covariate training data. This method posits that past behavior of the system is representative of its future behavior, so looking up states in a historical record that are closest to current observations is likely to yield a skillful forecast. Subsequent methodologies have also emphasized aspects of state space geometry, e.g., using the training data to approximate the evolution map through patched local linear models, often leveraging delay coordinates for state space reconstruction. 

Early approaches to coherent pattern extraction  include the proper orthogonal decomposition (POD) \cite{Kosambi43}, which is closely related to principal component analysis (PCA; introduced in the early 20th century by Pearson),  the Karhunen-Lo\`eve expansion, and empirical orthogonal function (EOF) analysis. Assuming that $ \mathcal Y $ is a Hilbert space, POD yields an expansion $ Y \approx Y_L = \sum_{j=1}^L z_j $, $ z_j = u_j \sigma_j \psi_j $.  Arranging the data into a matrix $\bm Y = ( y_0, \ldots, y_{N-1} )$, the $ \sigma_j $ are the singular values of $ \bm Y $ (in decreasing order), the $ u_j $ are the corresponding left singular vectors, called EOFs, and the $\psi_j$ are given by projections of $Y$ onto the EOFs, $\psi_j(\omega) = \langle u_j, Y(\omega) \rangle_{\mathcal Y}$. That is, the principal component $\psi_j : \Omega \to \mathbb R$ is a linear feature characterizing the unsupervised data $ \{ y_0, \ldots, y_{N-1} \} $. If the data is drawn from a probability measure $\mu$, as $N \to \infty$ the POD expansion is optimal in an $L^2(\mu)$ sense; that is,  $Y_L$ has minimal $L^2(\mu)$ error $ \lVert Y - Y_L \rVert_{L^2(\mu)}$ among all rank-$L$ approximations of $Y$. Effectively, from the perspective of POD, the important components of $Y$ are those capturing maximal variance.  

Despite many successes in challenging applications (e.g., turbulence), it has been recognized that POD may not reveal dynamically significant observables, offering limited predictability and physical insight. In recent years, there has been significant interest in techniques that address this shortcoming by modifying the linear map $ \bm Y $  to have an explicit dependence on the dynamics \cite{BroomheadKing86}, or replacing it by an evolution operator \cites{DellnitzJunge99,MezicBanaszuk99}. Either directly or indirectly, these methods make use of operator-theoretic ergodic theory, which we now discuss.

\section*{Operator-theoretic formulation}

The operator-theoretic formulation of dynamical systems theory shifts attention from the state-space perspective, and instead characterizes the dynamics through its action on linear spaces of observables. Denoting the vector space of $\mathcal Y$-valued functions on $\Omega$ by $\mathcal F$, for every time $t$ the dynamics has a natural induced action $U^t : \mathcal F \to \mathcal F$ given by composition with the flow map, $ U^t f = f \circ \Phi^t$. It then follows by definition that $U^t$ is a linear operator, i.e., $U^t( \alpha f + g ) = \alpha U^t f + U^ t g$ for all observables $f,g \in \mathcal F$ and every scalar $ \alpha \in \mathbb C$. The operator $U^t$ is known as a \emph{composition operator}, or \emph{Koopman operator} after  classical work of Bernard Koopman in the 1930s \cites{Koopman31}, which established that a general (potentially nonlinear) dynamical system can be characterized through intrinsically linear operators acting on spaces of observables. A related notion is that of the \emph{transfer operator}, $ P^t : \mathcal M \to \mathcal M $, which describes the action of the dynamics on a space of measures $\mathcal M$ via the pushforward map, $ P^t m := \Phi_{*}^t m = m \circ \Phi^{-t}  $.  In a number of cases, $\mathcal F$ and $\mathcal M $ are dual spaces to one another (e.g., continuous functions and Radon measures), in which case $U^t$ and $P^t$ are dual operators. 

If the space of observables under consideration is equipped with a Banach or Hilbert space structure, and the dynamics preserves that structure, the operator-theoretic formulation allows a broad range of tools from spectral theory and approximation theory for linear operators to be employed in the study of dynamical systems. For our purposes, a particularly advantageous aspect of this approach is that it is amenable to rigorous statistical approximation, which is one of our principal objectives. It should be kept in mind that the spaces of observables encountered in applications are generally infinite-dimensional, leading to behaviors with no counterparts in finite-dimensional linear algebra, such as unbounded operators and continuous spectrum. In fact, as we will see below, the presence of continuous spectrum is a hallmark of mixing (chaotic) dynamics.     

In this review, we restrict attention to the operator-theoretic description of  \emph{measure-preserving, ergodic dynamics}. By that, we mean that there is a probability measure $\mu$ on $\Omega$ such that (i) $\mu$ is invariant under the flow, i.e., $\Phi^t_*\mu = \mu$; and (ii) every measurable, $\Phi^t$-invariant set has either zero or full $\mu$-measure. We also assume that $\mu$ is a Borel measure with compact support $A \subseteq \Omega$; this set is necessarily $\Phi^t$-invariant. An example known to rigorously satisfy these properties is the Lorenz~63 (L63) system on $\Omega = \mathbb R^3$, which has a compactly supported, ergodic invariant measure supported on the famous ``butterfly'' fractal attractor; see Figure~\ref{figL63}. L63 exemplifies the fact that a smooth dynamical system may exhibit invariant measures with non-smooth supports. This behavior is ubiquitous in models of physical phenomena, which are formulated in terms of smooth differential equations, but whose long-term dynamics concentrate on lower-dimensional subsets of state space due to the presence of dissipation. Our methods should therefore not rely on the existence of a smooth structure for $A$.

In the setting of ergodic, measure-preserving dynamics on a metric space, two relevant structures that the dynamics may be required to preserve are continuity and $\mu$-measurability of observables. If the flow $ \Phi^t$ is continuous, then the Koopman operators act on the Banach space $\mathcal F = C(A,\mathcal Y)$ of continuous, $\mathcal Y$-valued functions on $A$, equipped with the uniform norm, by isometries, i.e., $\lVert U^t f \rVert_{\mathcal F} = \lVert f \rVert_{\mathcal F}$.  If $\Phi^t$ is $\mu$-measurable, then $U^t$ lifts to an operator on equivalence classes of $\mathcal Y $-valued functions in $L^p(\mu,\mathcal Y)$, $ 1 \leq p \leq \infty$, and acts again by isometries.  If $\mathcal Y$ is a Hilbert space (with inner product $\langle \cdot, \cdot \rangle_{\mathcal Y}$), the case $p = 2$ is special, since $  L^2(\mu,\mathcal Y)$ is a Hilbert space with inner product $ \langle f, g \rangle_{L^2(\mu,\mathcal Y)} = \int_\Omega \langle f(\omega),g(\omega) \rangle_{\mathcal Y} \, d\mu(\omega)$, on which $U^t$ acts as a unitary map, $U^{t*} = U^{-t}$. 

Clearly, the properties of approximation techniques for observables and evolution operators depend on the underlying space. For instance, $C(A,\mathcal Y)$ has a well-defined notion of pointwise evaluation at every $\omega \in \Omega$ by a continuous linear map $ \delta_\omega : C(A,\mathcal Y)\to \mathcal Y$, $ \delta_\omega f = f(\omega)$, which is useful for interpolation and forecasting,  but lacks an inner-product structure and associated orthogonal projections. On the other hand, $L^2(\mu)$ has inner-product structure, which is very useful theoretically as well as for numerical algorithms, but lacks the notion of pointwise evaluation. 

Letting $\mathcal F$ stand for any of the $C(A,\mathcal Y)$ or $ L^p(\mu,\mathcal Y)$ spaces,  the set $U = \{ U^t : \mathcal F \to \mathcal F \}_{t\in\mathbb R}$ forms a strongly continuous group under composition of operators. That is,  $U^t \circ U^s =  U^{t+s}$, $ U^{t,-1} = U^{-t} $, and $U^0 = \Id$, so that $U$ is a group, and for every $f \in \mathcal F$, $U^t f $ converges to $f $ in the norm of $\mathcal F$ as $ t \to 0 $. A central notion in such evolution groups is that of the \emph{generator}, defined by the $\mathcal F$-norm limit $ V f = \lim_{t\to 0} ( U^t f - f ) / t $ for all $ f \in \mathcal F$ for which the limit exists. It can be shown that the domain $D(V)$ of all such $f$ is a dense subspace of $\mathcal F$, and $V : D(V) \to \mathcal F$ is a closed, unbounded operator. Intuitively,  $V$ can be thought as a directional derivative of observables along the dynamics. For example, if $\mathcal Y = \mathbb C$, $A$ is a $C^1$ manifold, and the flow $ \Phi^t : A \to A$ is generated by a continuous vector field $ \vec V : A \to TA$, the generator of the Koopman group on $C(A)$ has as its domain the space $C^1(A) \subset C(A)$ of continuously differentiable, complex-valued functions, and $V f = \vec V \cdot \nabla f$ for $ f \in C^1(A)$. A strongly continuous evolution group is completely characterized by its generator, as any two such groups with the same generator are identical. 

The generator acquires additional properties in the setting of unitary evolution groups on $H = L^2(\mu,\mathcal Y)$, where it is skew-adjoint, $V^* = - V$. Note that the skew-adjointness of $V$ holds for more general measure-preserving dynamics than Hamiltonian systems, whose generator is skew-adjoint with respect to Lebesgue measure. By the spectral theorem for skew-adjoint operators, there exists a unique projection-valued measure $E : \mathcal B(\mathbb R) \to B(H)$, giving the generator and Koopman operator as the spectral integrals
\begin{displaymath}
    V = \int_{\mathbb R} i \alpha \, dE(\alpha), \quad U^t = e^{tV} = \int_{\mathbb R} e^{i\alpha t}\, dE(\alpha).
\end{displaymath}
Here, $\mathcal B(\mathbb R)$ is the Borel $\sigma$-algebra on the real line, and $B(H)$ the space of bounded operators on $H$.  Intuitively,  $E$ can be thought of as an operator analog of a complex-valued spectral measure in Fourier analysis, with $\mathbb R$ playing the role of frequency space.  That is, given $ f \in H $, the $\mathbb C$-valued Borel measure $E_f(S) = \langle f, E(S) f \rangle_{H} $ is precisely the Fourier spectral measure associated with the time-autocorrelation function $ C_f(t) =  \langle f, U^t f \rangle_H$. The latter, admits the Fourier representation $C_f(t) = \int_{\mathbb R} e^{i\alpha t} \, dE_f(\alpha)$.

The Hilbert space $H$ admits a $U^t$-invariant splitting $H = H_a \oplus H_c$ into orthogonal subspaces $H_a$ and $H_c$ associated with the point and continuous components of $E$, respectively. In particular, $E$ has a unique decomposition $E= E_a + E_c $ with $ H_a = \ran E_a( \mathbb R )$ and $ H_c = \ran E_c(\mathbb R) $, where $E_a$ is a purely atomic spectral measure, and $E_c$ is a spectral measure with no atoms. The atoms of $E_a$ (i.e., the singletons $ \{ \alpha_j \} $ with $E_a(\{ \alpha_j \} ) \neq 0 $) correspond to \emph{eigenfrequencies} of the generator, for which the eigenvalue equation $ V z_j = i \alpha z_j $  has a nonzero solution $z_j \in H_a$.  Under ergodic dynamics, every eigenspace of $V$ is one-dimensional, so that if $z_j$ is normalized to unit $L^2(\mu)$ norm, $E(\{ \alpha_j \} ) f = \langle z_j,f \rangle_{L^2(\mu)} z_j $. Every such $z_j$  is an eigenfunction of the Koopman operator $U^t$ at eigenvalue $e^{i\alpha_j t}$, and $ \{ z_j \} $ is an orthonormal basis of $H_a$. Thus, every $ f \in H_a $ has the quasiperiodic evolution $U^t f = \sum_j e^{i\alpha_j t} \langle z_j, f \rangle_{L^2(\mu)} z_j$, and the autocorrelation $C_f(t)$ is also quasiperiodic. While $H_a$ always contains constant eigenfunctions with zero frequency, it might not have any non-constant elements. In that case, the dynamics is said to be \emph{weak-mixing}. In contrast to the quasiperiodic evolution of observables in $H_a$, observables in the continuous spectrum subspace exhibit a loss of correlation characteristic of mixing (chaotic) dynamics. Specifically, for every $ f \in H_c$ the time-averaged autocorrelation function $ \bar C_f(t) = \int_{0}^t \lvert C_f(s) \rvert \, ds / t $ tends to 0 as $ \lvert t \rvert \to \infty$, as do cross-correlation functions $ \langle g, U^t f \rangle_{L^2(\mu)}$ between observables in $H_c$ and arbitrary observables in $L^2(\mu)$.

\section*{Data-driven forecasting}

Based on the concepts introduced above, one can formulate statistical forecasting in Problem~\ref{probForecast} as the task of constructing a function $Z_t : \mathcal X \to \mathcal Y$ on covariate space $\mathcal X$, such that $ Z_t \circ X $ optimally approximates $U^t Y $ among all functions in a suitable class. We set $\mathcal Y = \mathbb C$, so the response variable is scalar-valued, and consider the Koopman operator on $L^2(\mu)$, so we have access to orthogonal projections. We also assume for now that the covariate function $X $ is injective, so $ \hat Y_t := Z_t \circ X $ should be able to approximate $U^t Y $ to arbitrarily high precision in $L^2(\mu) $ norm. Indeed, let $ \{ u_0, u_1, \ldots \} $ be an orthonormal basis of $L^2(\nu)$, where $\nu = X_* \mu $ is the pushforward of the invariant measure onto $\mathcal X$. Then, $\{ \phi_0, \phi_1, \ldots \}$ with $ \phi_j = u_j \circ X $ is an orthonormal basis of $L^2(\mu)$. Given this basis, and because $U^t$ is bounded, we have $ U^tY = \lim_{L\to\infty} U^t_L Y$, where the partial sum $U^t_LY := \sum_{j=0}^{L-1}  \langle U^t Y,\phi_j \rangle_{L^2(\mu)} \phi_j$ converges in $L^2(\mu)$ norm. Here, $ U^t_L $ is a finite-rank map on $L^2(\mu)$ with range $ \spn\{ \phi_0, \ldots, \phi_{L-1} \}$, represented by  an $L\times L$ matrix $\bm U(t)$ with elements $U_{ij}(t) = \langle \phi_i, U^t \phi_j \rangle_{L^2(\mu)}$. Defining $ \vec y = ( \hat y_0, \ldots, \hat y_{L-1} )^\top$, $\hat y_j = \langle \phi_j, U^t Y \rangle_{L^2(\mu)}$, and $ ( \hat z_0(t), \ldots, \hat z_{L-1}(t) )^\top = \bm U(t) \vec y$, we have $ U^t_L Y =  \sum_{j=0}^{L-1} \hat{z}_j(t) \phi_j $. Since $\phi_j = u_j \circ X$, this leads to the estimator $ \hat Z_{t,L} \in L^2(\nu) $, with $\hat Z_{t,L} = \sum_{j=0}^{L-1} \hat z_j(t) u_j $. 

The approach outlined above tentatively provides a consistent forecasting framework. Yet, while in principle appealing, it has three major shortcomings: (i) Apart from special cases, the invariant measure and an orthonormal basis of $L^2(\mu)$ are not known. In particular, orthogonal functions with respect to an ambient measure on $\Omega$ (e.g., Lebesgue-orthogonal polynomials) will not suffice, since there are no guarantees that such functions form a Schauder basis of $L^2(\mu)$, let alone be orthonormal. Even with a basis, we  cannot evaluate $U^t$ on its elements without knowing $ \Phi^t$. (ii) Pointwise evaluation on $L^2(\mu)$ is not defined, making $\hat Z_{t,L}$ inadequate in practice, even if the coefficients $\hat z_j(t)$ are known. (iii) The covariate map $X$ is oftentimes non-invertible, and thus the $ \phi_j $ span a strict subspace of $L^2(\mu)$. We now describe methods to overcome these obstacles using learning theory.

\subsection*{Sampling measures and ergodicity}

The dynamical trajectory  $ \{ \omega_0, \ldots, \omega_{N-1} \} $ in state space underlying the training data is the support of a discrete \emph{sampling measure}  $\mu_N := \sum_{n=0}^{N-1} \delta_{\omega_n} / N$.  A key consequence of ergodicity is that for Lebesgue-a.e.\ sampling interval $\Delta t$ and $\mu$-a.e.\ starting point $\omega_0 \in \Omega$, as $N \to \infty$, the sampling measures $\mu_N$ weak-converge to the invariant measure $\mu$; that is,
\begin{equation}
    \label{eqWeakConv}
    \lim_{N\to \infty} \int_\Omega f \, d\mu_N =  \int_\Omega f \, d\mu, \quad \forall f \in C(\Omega).
\end{equation}
Since integrals against $\mu_N$ are time averages on dynamical trajectories, $ \int_{\Omega} f \, d\mu_N = \sum_{n=0}^{N-1} f(\omega_n) / N $,  ergodicity provides an empirical means of accessing the statistics of the invariant measure. In fact,  many systems encountered in  applications possess so-called \emph{physical measures}, where~\eqref{eqWeakConv} holds for $\omega_0$ in a ``larger'' set of positive measure with respect to an ambient measure (e.g., Lebesgue measure)  from which experimental initial conditions are drawn. Hereafter, we will let $M$ be a compact subset of $\Omega$, which is forward-invariant under the dynamics (i.e., $\Phi^t(M) \subseteq M$ for all $ t \geq 0 $), and thus necessarily contains $A$.  For example, in dissipative dynamical systems such as L63, $M$ can be chosen as a compact absorbing ball. 

\subsection*{Shift operators}

Ergodicity suggests that appropriate data-driven analogs of are the $L^2(\mu_N)$ spaces induced by the the sampling measures $ \mu_N$. For a given $N$, $L^2(\mu_N)$ consists of equivalence classes of measurable functions $f : \Omega \to \mathbb C$ having common values at the sampled states $\omega_n$, and the inner product of two elements $f,g \in L^2(\mu_N)$ is given by an empirical time-correlation, $ \langle f, g \rangle_{\mu_N} = \int_\Omega f^* g \, d\mu_N = \sum_{n=0}^{N-1} f^*(\omega_n) g(\omega_n) / N $. Moreover, if the $\omega_n$ are distinct (as we will assume for simplicity of exposition), $L^2(\mu_N)$ has dimension $N$, and is isomorphic as a Hilbert space to $\mathbb C^N$ equipped with a normalized dot product. Given that, we can represent every $ f \in L^2(\mu_N)$ by a column vector $ \vec f = ( f(\omega_0), \ldots, f(\omega_{N-1}))^\top \in \mathbb C^N $, and every linear map $A : L^2(\mu_N) \to L^2(\mu_N)$ by an $N\times N$ matrix $\bm A$, so that $ \vec g = \bm A \vec f$ is the column vector representing $g = A f$. The elements of  $ \vec f  $ can also be understood as expansion coefficients in the standard basis $ \{ e_{0,N}, \ldots, e_{N-1,N} \} $ of $L^2(\mu_N) $, where $ e_{j,N}(\omega_n) = N^{1/2}\delta_{jn} $; that is, $f(\omega_n) = \langle e_{n,N},f \rangle_{L^2(\mu_N)}$. Similarly, the elements of $\bm A $ correspond to the operator matrix elements $A_{ij} = \langle e_{i,N}, A e_{j,N} \rangle_{L^2(\mu_N)}$. 



Next, we would like to define a Koopman operator on $L^2(\mu_N)$, but this space does not admit a such an operator as a composition map induced by the dynamical flow $ \Phi^t$ on $\Omega$. This is because $ \Phi^t $ does not preserve null sets with respect to $\mu_N$, and thus does not lead to a well-defined composition map on equivalence classes of functions in $L^2(\mu_N)$. Nevertheless, on $L^2(\mu_N)$ there is an analogous construct to the Koopman operator on $L^2(\mu)$, namely the \emph{shift operator}, $U^q_{N} : L^2(\mu_N) \to L^2(\mu_N)$, $ q \in \mathbb Z $, defined as
\begin{displaymath}
    U^q_N f(\omega_n) = 
    \begin{cases}
        f(\omega_{n+q}), & 0 \leq n+q \leq N -1, \\
        0, & \text{otherwise}.
    \end{cases}
\end{displaymath}
Even though $U^q_N$ is not a composition map, intuitively it should have a connection with the Koopman operator $U^{q \, \Delta t}$. One could consider, for instance, the matrix representation $ \tilde{ \bm U}_N(q) = [ \langle e_{i,N}, U^q_N e_{j,N} \rangle_{L^2(\mu_N)} ]$ in the standard basis, and attempt to connect it with a matrix representation of $U^{q\,\Delta t}$ in an orthonormal basis of $L^2(\mu)$. However, the issue with this approach is that the  $e_{j,N}$ do not have $N\to \infty $ limits in $L^2(\mu)$, meaning that there is no suitable notion of $N\to \infty$ convergence of the matrix elements of $U^q_N$ in the standard basis.  In response,  we will construct a representation of the shift operator in a different orthonormal basis with a well-defined $N\to\infty$ limit. The main tools that we will use  are \emph{kernel integral operators}, which we now describe. 
\subsection*{Kernel integral operators}

In the present context, a \emph{kernel function} will be a real-valued, continuous function $k : \Omega \times \Omega \to \mathbb R$ with the property that there exists a strictly positive, continuous function $ d : \Omega \to \mathbb R$ such that 
\begin{equation}
    \label{eqBalance}
    d( \omega ) k( \omega, \omega' ) = d( \omega' )k( \omega', \omega ), \quad \forall \omega,\omega' \in \Omega.    
\end{equation}
Notice the similarity between~\eqref{eqBalance} and the detailed balance relation in reversible Markov chains. Let now $\rho$ be any Borel probability measure with compact support $S \subseteq M$ included in the forward-invariant set $M$. It follows by continuity of $k$ and compactness of $S$ that the integral operator $K_\rho : L^2(\rho) \to C(M)$, 
\begin{equation}
    \label{eqKOp}
    K_\rho f = \int_\Omega k( \cdot, \omega ) f(\omega)\, d\rho(\omega),
\end{equation}
is well-defined as a bounded operator mapping elements of $L^2(\rho)$ into continuous functions on $M$. Using $\iota_\rho : C(M) \to L^2(\rho)$ to denote the canonical inclusion map, we consider two additional integral operators, $G_\rho : L^2(\rho) \to L^2(\rho)$ and $\tilde G_\rho : C(M) \to C(M)$, with $G_\rho = \iota_\rho K_\rho $ and $ \tilde G_\rho = K_\rho \iota_\rho$, respectively. 

The operators $G_\rho$ and $\tilde G_\rho$ are compact operators acting with the same integral formula as $K_\rho$ in~\eqref{eqKOp}, but their codomains and domains, respectively, are different. Nevertheless, their nonzero eigenvalues coincide, and $\phi \in L^2(\rho) $ is an eigenfunction of $G_\rho $ corresponding to a nonzero eigenvalue $\lambda$ if and only if $ \varphi \in C(M) $ with $ \varphi = K_\rho \phi / \lambda$ is an eigenfunction of $\tilde G_\rho $ at the same eigenvalue. In effect, $\phi \mapsto \varphi $ ``interpolates'' the $L^2(\rho)$ element $\phi$ (defined only up to $\rho$-null sets) to the continuous, everywhere-defined function $\varphi$. It can be verified that if~\eqref{eqBalance} holds, $G_\rho$ is a trace-class operator with real eigenvalues,  $ \lvert \lambda_0 \rvert \geq \lvert \lambda_1 \rvert \geq \cdots \searrow 0^+$. Moreover, there exists a Riesz basis $ \{ \phi_0, \phi_1, \ldots, \} $ of $L^2(\rho)$ and a corresponding dual basis $ \{ \phi'_0, \phi'_1, \ldots \} $ with $ \langle \phi'_i, \phi_j \rangle_{L^2(\rho)} = \delta_{ij} $, such that $ G_\rho \phi_j = \lambda_j \phi_j $ and $ G_\rho^* \phi'_j = \lambda_j \phi'_j $. We say that the kernel $k$ is \emph{$L^2(\rho)$-universal} if $G_\rho$ has no zero eigenvalues; this is equivalent to  $ \ran G_\rho$ being dense in $L^2(\rho)$. Moreover, $k$ is said to be \emph{$L^2(\rho)$-Markov} if $G_\rho$ is a Markov operator, i.e., $G_\rho \geq 0 $, $ G_\rho f \geq 0 $ if $ f \geq  0 $, and $ G 1 = 1 $. 

Observe now that the operators $G_{\mu_N}$ associated with the sampling measures $\mu_N$, henceforth abbreviated by $G_N$, are represented by $N\times N$ kernel matrices $\bm G_N = [ \langle e_{i,N}, G_N e_{j,N} \rangle_{L^2(\mu_N)}] = [ k(\omega_i,\omega_j) ]  $ in the standard basis of $L^2(\mu_N)$. Further, if $k$ is a pullback kernel from covariate space, i.e., $k(\omega,\omega') = \kappa(X(\omega),X(\omega'))$ for  $ \kappa : \mathcal X \times \mathcal X \to \mathbb R$, then $\bm G_N = [ \kappa(x_i,x_j) ]$ is empirically accessible from the training data.  Popular kernels in applications include the covariance kernel $ \kappa(x,x') = \langle x, x' \rangle_{\mathcal X}$ on an inner-product space and the radial Gaussian kernel $ \kappa(x,x') = e^{-\lVert x- x' \rVert^2_{\mathcal  X}/\epsilon}$ \cite{Genton01}. It is also common to employ Markov kernels constructed by normalization of symmetric kernels \cites{CoifmanLafon06,BerryHarlim16}. We will use $k_N$ to denote kernels with data-dependent normalizations. 
    
\begin{figure*}
    \centering
    \includegraphics[width=.84\linewidth]{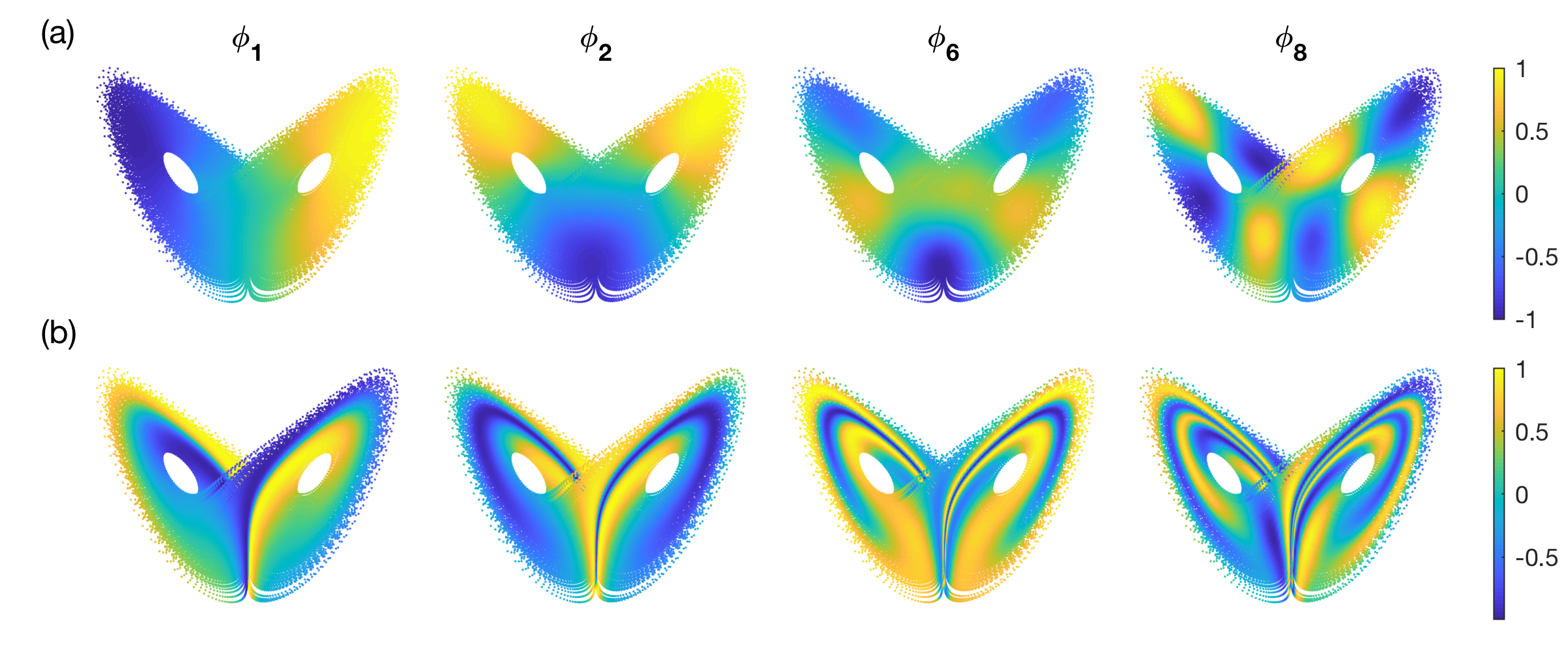}\\
    \includegraphics[width=.74\linewidth]{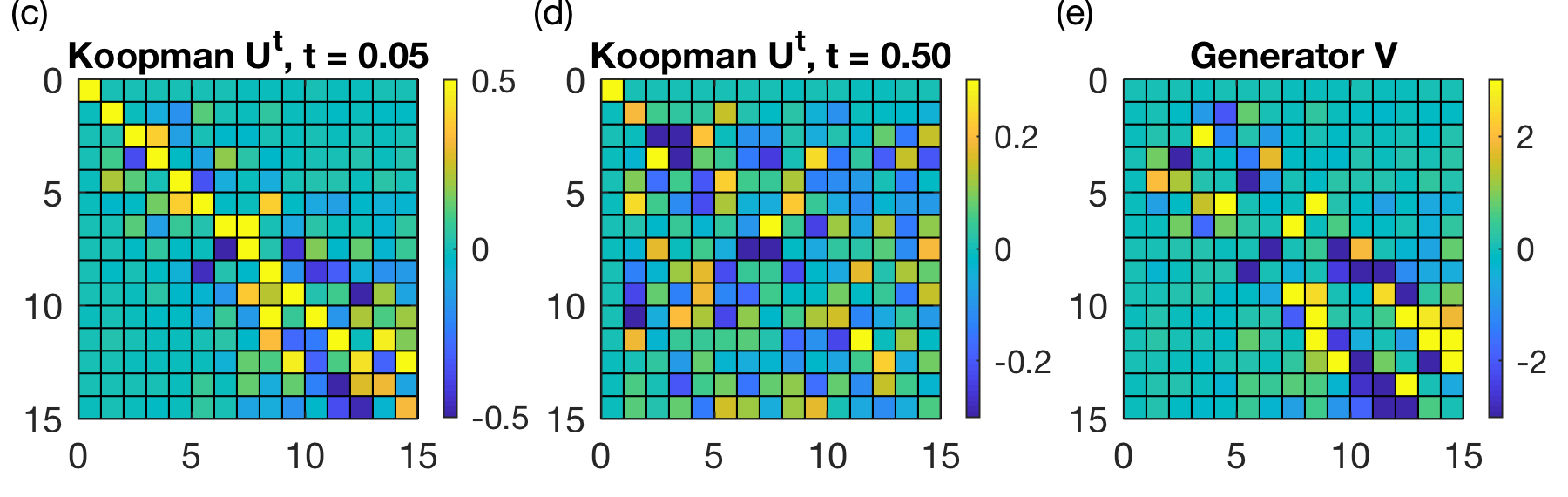}
    \small

    \begin{tikzpicture} 
    \node[draw,label={[label distance=.01cm]174:(f)}]{\begin{tikzcd}[column sep=2.5em,row sep=0em]
            & \mathcal{H}_N \arrow{r}{\iota_N} & L^2(\mu_N) \arrow{r}{U^{q*}_{N}}   & L^2(\mu_N) \arrow{r}{\Pi_L} &\textup{span}(\{\phi_{i,N}\}_{i=1}^L) \arrow{r}{\mathbb{E}_{(\cdot)}Y} & \mathcal{Y} \arrow[dashrightarrow]{dd}{\textup{error}} \\
       M\subseteq\Omega \arrow{ru}{\Psi_N}\arrow{rd}{\Psi}  & & & & & \\
       & \mathcal{H}(M) \arrow{r}{\iota} & L^2(\mu) \arrow{r}{U^{t*}}  & L^2(\mu)  \arrow{rr}{\mathbb{E}_{(\cdot)} Y} & & \mathcal{Y} \\
       \omega \arrow{uu}{\in} \arrow[mapsto]{r} &\Psi(\omega)  \arrow[mapsto]{rr} & & U^{t*}\Psi(\omega)  \arrow[mapsto]{rr} & & \mathbb{E}_{\Psi(\omega)}U^{t} Y 
    \end{tikzcd}};
\end{tikzpicture}\vspace{5pt}
\begin{tikzpicture}
\node[draw,label={[label distance=.01cm]174:(g)}]{
    \begin{tikzcd}[column sep=2.73em,row sep=0em]
       	& L^2(\mu_N) \arrow{r}{U^q_N} & L^2(\mu_N) \arrow{r}{\Pi_X}   & L^2_X(\mu_N) \arrow{r}{\mathcal{N}_N} & \mathcal{H}_N \arrow{r}{\iota} & L^2_X(\mu) \arrow[dashrightarrow]{dd}{\textup{error}} \\
       C(M) \arrow{ru}{\iota_N}\arrow{rd}{\iota}  & & & & & \\
       & L^2(\mu) \arrow{r}{U^{t}} & L^2(\mu) \arrow{rrr}{\Pi_X}  & & & L^2_X(\mu) \\
       Y \arrow{uu}{\in} \arrow[mapsto]{rr} & & U^tY \arrow[mapsto]{rrr}   & & & Z_{t} \circ X = \mathbb{E}(U^{t}Y \mid X) 
    \end{tikzcd}};
\end{tikzpicture}

\caption{\label{figL63}  Panel~(a) shows eigenfunctions $\phi_{j,N}$ of $G_N$ for a dataset sampled near the L63 attractor. Panel~(b) shows the action of the shift operator $U^q_N$ on the $\phi_{j,N}$ from (a) for $q=50$ steps, approximating the Koopman operator $U^{q\,\Delta t} $. Panels~(c, d) show the matrix elements $\langle \phi_{i,N}, U^q_N \phi_{j,N} \rangle_{\mu_N}$ of the shift operator for $ q = 5$ and 50. The mixing dynamics is evident in the larger far-from-diagonal components in $q=50$ vs.\ $q=5$. Panel~(e) shows the matrix representation of a finite-difference approximation of the generator $V$, which is skew-symmetric.  Panels~(f, g) summarize the diffusion forecast (DF) and kernel analog forecast (KAF) for lead time $ t = q \, \Delta t$. In each diagram, the data-driven finite dimensional approximation (top row) converges to the true forecast (middle row). DF maps an initial state $\omega \in M \subseteq \Omega$ to the future expectation of an observable $\mathbb{E}_{\Psi(\omega)}U^t Y = \mathbb E_{U^{t*}\Psi(\omega)} Y $ and KAF maps a response function $Y \in C(M)$ to the conditional expectation $\mathbb E(U^tY \mid X)$.} 
    
\end{figure*}
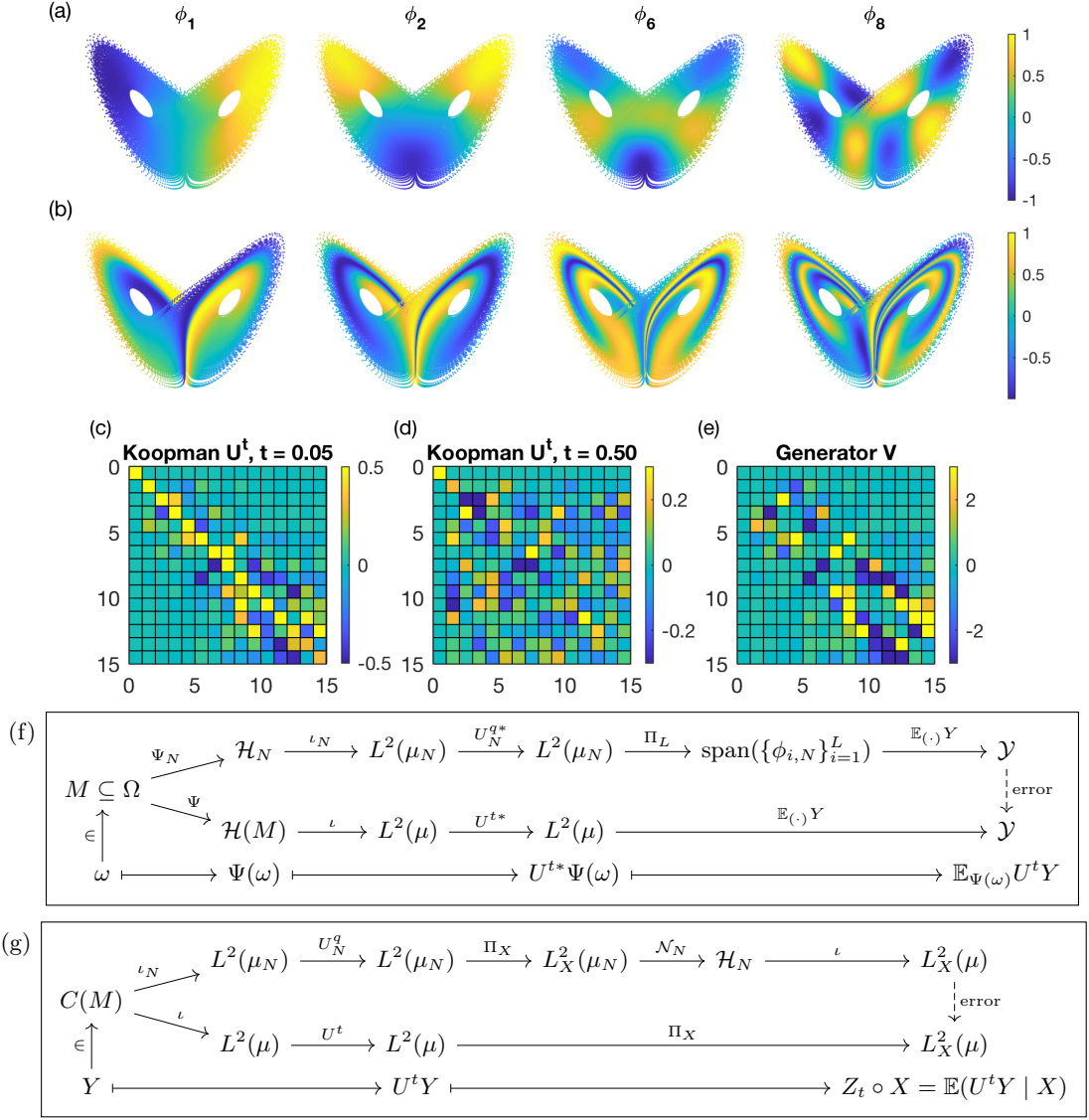

A widely used strategy for learning with integral operators  \cite{VonLuxburgEtAl08} is to construct families of kernels $k_N$ converging in $C(M\times M) $ norm to $k$. This implies that for every nonzero eigenvalue $\lambda_j$ of $G \equiv G_\mu $, the sequence of eigenvalues $\lambda_{j,N}$ of $G_N $ satisfies $\lim_{N\to\infty \lambda_{j,N}} = \lambda_j$. Moreover, there exists a sequence of eigenfunctions $\phi_{j,N} \in L^2(\mu_N) $ corresponding to $\lambda_{j,N}$, whose continuous representatives, $\varphi_{j,N} = K_N \phi_{j,N} / \lambda_{j,N}$, converge in $C(M)$ to $ \varphi_j = K \phi_j / \lambda_j $, where $ \phi_j \in L^2(\mu) $ is any eigenfunction of $G$ at eigenvalue $\lambda_j$. In effect, we use $C(M)$  as a ``bridge'' to establish spectral convergence of the operators $G_N$, which act on different spaces. Note that  $(\lambda_{j,N},\varphi_{j,N})$ does not converge uniformly with respect to $ j$, and for a fixed $N$, eigenvalues/eigenfunctions at larger $j$ exhibit larger deviations from their $N\to\infty$ limits. Under measure-preserving, ergodic dynamics, convergence occurs for $ \mu $-a.e.\ starting state $ \omega_0 \in M$, and $\omega_0 $ in a set of positive ambient measure if $ \mu $ is physical. In particular, the training states $\omega_n$ need not lie on $A$. See Figure~\ref{figL63} for eigenfunctions of $G_N$ computed from data sampled near the L63 attractor. When the invariant measure $ \mu $ has a smooth density with respect to local coordinates on $ \Omega $, results on spectral convergence of graph Laplacians to manifold Laplacians \cites{TrillosSlepcev18,TrillosEtAl19} could be employed to provide a more precise characterization of the spectral properties of $G_N$ for suitable choices of kernel.

\subsection*{Diffusion forecasting}

We now have the ingredients to build a concrete statistical forecasting scheme based on data-driven approximations of the  Koopman operator. In particular, note that if $\phi'_{i,N},\phi_{j,N} $ are biorthogonal eigenfunctions of $G^*_N$ and $G_N$, respectively, at nonzero eigenvalues,  we can evaluate the matrix element $ U_{N,ij}(q): = \langle \phi'_{i,N}, U^q_N \phi_{j,N} \rangle_{L^2(\mu_N)}$ of the shift operator using the continuous representatives $\varphi'_{i,N},\varphi_{j,N} $,  
\begin{align*}
    U_{N,ij}(q) &= \frac{1}{N}\sum_{n=0}^{N-1-q} \phi'_{i,N}(\omega_n) \phi_{j,N}(\omega_{n+q})\\
    &= \frac{N - q}{N}\int_\Omega \varphi'_{i,N} U^{q\,\Delta t} \varphi_{j,N} \, d\mu_{N-q},
\end{align*}
where $U^{q\,\Delta t} $ is the Koopman operator on $C(M)$. Therefore, if the corresponding eigenvalues $\lambda_i,\lambda_j$ of $G$ are nonzero, by the weak convergence of the sampling measures in~\eqref{eqWeakConv} and uniform convergence of the eigenfunctions, as $N\to\infty$, $ U_{ij,N}(q)$  converges to the matrix element $U_{ij}(q\, \Delta t) = \langle \phi_i, U^{q\,\Delta t} \phi_j \rangle_{L^2(\mu)} $  of the Koopman operator on $L^2(\mu)$. This convergence is not uniform with respect to $i,j$, but if we fix a parameter $ L \in \mathbb N$ (which can be thought of as spectral resolution) such that $\lambda_{L-1} \neq 0$, we can obtain a statistically consistent approximation of $L\times L$ Koopman operator matrices, $\bm U(q\,\Delta t) = [ U_{ij}(q\,\Delta t)]$, by shift operator matrices, $\bm U_N(q) = [U_{N,ij}(q)]$, with $i,j \in \{ 0, \ldots, L-1\}$. Checkerboard plots of $\bm U_N(q)$ for the L63 system are displayed in Figure~\ref{figL63}.    

This method for approximating matrix elements of Koopman operators was proposed in a technique called \emph{diffusion forecasting} (named after the diffusion kernels employed) \cite{BerryEtAl15}. Assuming that the response $Y$ is continuous and  by spectral convergence of $G_N$, for every $j \in \mathbb N_0$ such that $\lambda_j > 0$, the inner products  $\hat Y_{j,N} = \langle \phi'_{j,N}, Y \rangle_{\mu_N} $  converge, as $N\to\infty$, to $ \hat Y_j = \langle \phi'_j, Y \rangle_{L^2(\mu)}$. This implies that for any $L \in \mathbb N$ such that $\lambda_{L-1} > 0 $,  $ \sum_{j=0}^{L-1} \hat Y_{j,N} \varphi_{j,N} $ converges in $C(M)$  to the continuous representative of $ \Pi_L Y$, where $\Pi_L $ is the orthogonal projection on $L^2(\mu)$ mapping into $ \spn\{\phi_0, \ldots, \phi_{L-1} \} $. Suppose now that $\varrho_N$ is a sequence of continuous functions converging uniformly to $\varrho \in C(M)$, such that $\varrho_N$ are probability densities with respect to $\mu_N$ (i.e., $\varrho_N \geq 0 $ and $\lVert \varrho_N \rVert_{L^1(\mu_N)} = 1$). By similar arguments as for $Y$, as $N\to\infty$, the continuous function $ \sum_{j=0}^{L-1} \hat \varrho_{j,N} \varphi_{j,N}$ with $ \hat \varrho_{j,N} = \langle \varphi'_{j,N}, \varrho_N\rangle_{L^2(\mu_N)}$ converges to $ \Pi_L \varrho$ in $L^2(\mu)$. Putting these facts together, and setting $ \vec \varrho_N = ( \hat \varrho_{0,N}, \ldots, \hat \varrho_{L-1,N} )^\top$ and $ \vec Y_N = ( \hat Y_{0,N}, \ldots, \hat Y_{L-1,N})^\top$, we conclude that
\begin{equation}
    \label{eqDF}
    \vec \varrho_N^\top \bm U_N(q) \vec Y_N \xrightarrow{N\to\infty} \langle \Pi_L \varrho, \Pi_L U^{q\,\Delta t} Y \rangle_{L^2(\mu)}. 
\end{equation}
Here, the left-hand side is given by matrix--vector products obtained from the data,  and the right-hand side is equal to the expectation of  $\Pi_L U^{q\,\Delta t} Y$ with respect to the probability measure $ \rho $ with density $ d\rho/d\mu = \varrho$; i.e., $ \langle \Pi_L \varrho, \Pi_L U^{q\,\Delta t} Y \rangle_{L^2(\mu)} = \mathbb E_\rho( \Pi_L U^{q\,\Delta t} Y )$, where $ \mathbb E_\rho(\cdot) := \int_\Omega (\cdot)\,  d\rho $.

What about the dependence of the forecast on $L$? As $L$ increases, $\Pi_L$ converges strongly to the orthogonal  projection $\Pi_G : L^2(\mu) \to L^2(\mu)$ onto the closure of the range of $G$. Thus, if the kernel $k$ is $L^2(\mu)$-universal (i.e., $\overline{\ran G} = L^2(\mu)$), $\Pi_G = \Id$, and under the iterated limit of $L\to \infty$ after $N \to \infty $ the left-hand side of~\eqref{eqDF} converges to $\mathbb E_\rho U^{q\,\Delta t} Y$. In summary, implemented with an $L^2(\mu)$-universal kernel, diffusion forecasting consistently approximates the expected value of the time-evolution of any continuous observable with respect to any probability measure with continuous density relative to $\mu$. An example of an $L^2(\mu)$-universal kernel is the pullback of a radial Gaussian kernel on  $\mathcal X = \mathbb R^m$. In contrast, the covariance kernel is not $L^2(\mu)$-universal, as in this case the rank of $G$ is bounded by $m$. This illustrates that forecasting in the POD basis may be subject to intrinsic limitations, even with full observations.   


\subsection*{Kernel analog forecasting}

While providing a flexible framework for approximating expectation values of observables under measure-preserving, ergodic dynamics, diffusion forecasting does not directly address the problem of constructing a concrete forecast function, i.e., a function $Z_t : \mathcal X \to \mathbb C$ approximating $U^t Y$ as stated in Problem~\ref{probForecast}. One way of defining such a function is to let $\kappa_N$ be a $L^2(\nu_N)$-Markov kernel on $\mathcal X $ for $\nu_N = X_* \mu_N$, and to consider the ``feature map'' $ \Psi_N : \mathcal X \to C(M) $ mapping each point $ x \in \mathcal X$ in covariate space to the kernel section $ \Psi_N(x) = \kappa_N(x, X(\cdot))$. Then, $ \Psi_N(x)$ is a continuous probability density with respect to $\mu_N$, and we can use diffusion forecasting to define $Z_{q\,\Delta t}(x) = \overrightarrow{\Psi_N(x)}^\top \bm U_N(q) \vec Y_N$ with notation as in~\eqref{eqDF}. 

While this approach has a well-defined $N\to\infty$ limit, it does not provide optimality guarantees, particularly in situations where $X$ is non-injective. Indeed, the $L^2(\mu)$-optimal approximation to $U^tY$ of the form $Z_t \circ X$ is given by the \emph{conditional expectation} $\mathbb E(U^tY \mid X)$. In the present, $L^2 $, setting we have $ \mathbb E(U^tY \mid X) = \Pi_X U^tY $, where $ \Pi_X $ is the orthogonal projection into $L^2_X(\mu) := \{ f \in L^2(\mu): f = g \circ X \} $. That is, the conditional expectation minimizes the error $ \lVert f - U^t Y \rVert^2_{L^2(\mu)}$ among all pullbacks $f \in L^2_X(\mu) $ from covariate space. Even though $\mathbb E(U^tY \mid X = x)$ can be expressed as an expectation with respect to a conditional probability measure $\mu(\cdot \mid x ) $ on $\Omega$, that measure will generally not have an  $L^2(\mu)$ density, and there is no map $\Psi : \mathcal X \to C(M)$  such that $ \langle \Psi(x), U^t Y \rangle_{L^2(\mu)} $ equals $\mathbb E(U^tY \mid X = x)$.    

To construct a consistent estimator of the conditional expectation, we require that $k$ is a pullback of a kernel $\kappa : \mathcal X \times \mathcal X \to \mathbb R$ on covariate space which is (i) symmetric, $\kappa(x,x') = \kappa(x', x) $ for all $x,x' \in \mathcal X$ (so \eqref{eqBalance} holds); (ii) strictly positive;  and (iii) \emph{strictly positive-definite}. The latter means that for any sequence $x_0, \ldots, x_{n-1}$ of distinct points in $\mathcal X$ the  matrix $[\kappa(x_i,x_j)]$ is strictly positive. These properties imply that there exists a Hilbert space $\mathcal H$ of complex-valued functions on $\Omega$, such that (i) for every $\omega\in \Omega$, the kernel sections $k_\omega = k(\omega,\cdot)$ lie in $\mathcal H$; (ii) the evaluation functional $\delta_\omega : \mathcal H \to \mathbb C$ is bounded and satisfies $\delta_\omega f = \langle k_\omega, f \rangle_{\mathcal H}$; (iii) every $f \in \mathcal H$ has the form $ f = g \circ X $ for a continuous function $ g : \mathcal X \to \mathbb C$; and (iv) $\iota_\mu \mathcal H $ lies dense in $L^2_X(\mu)$. 

A Hilbert space of functions satisfying (i) and (ii) above is known as a \emph{reproducing kernel Hilbert space (RKHS)}, and the associated kernel $k$ is known as a \emph{reproducing kernel}. RKHSs have many useful properties for statistical learning \cite{CuckerSmale01}, not least because they combine the Hilbert space structure of $L^2$ spaces with pointwise evaluation in spaces of continuous functions.  The density of $\mathcal H$ in $L^2_X(\mu) $ is a consequence of the strict positive-definiteness of $\kappa$. In particular, because the conditional expectation $\mathbb E(U^tY \mid X)$ lies in $L^2_X(\mu)$, it can be approximated by elements of $\mathcal H$ to arbitrarily high precision in $L^2(\mu)$ norm, and every such approximation will be a pullback $ \hat Y_t = Z_t \circ X$ of a continuous function $Z_t $ that can be evaluated at arbitrary covariate values.  



We now describe a data-driven technique for constructing such a prediction function, which we refer to as \emph{kernel analog forecasting (KAF)} \cite{AlexanderGiannakis20}. Mathematically, KAF closely related to kernel principal component regression. To build the KAF estimator, we work again with integral operators as in~\eqref{eqKOp}, with the difference that now $K_\rho : L^2(\rho) \to \mathcal H(M)$ takes values in the restriction of $\mathcal H$ to the forward-invariant set $M$, denoted $\mathcal H(M)$.  One can show that the adjoint $K^*_\rho : \mathcal H(M) \to L^2(\rho)$ coincides with the inclusion map $ \iota_\rho $ on continuous functions, so that $K^*_\rho $ maps $ f \in \mathcal H(M) \subset C(M)$ to its corresponding $L^2(\rho)$ equivalence class. As a result, the integral operator $G_\rho : L^2(\rho) \to L^2(\rho)$ takes the form $ G_\rho = K^*_\rho K_\rho$, becoming a self-adjoint, positive-definite, compact operator with eigenvalues $ \lambda_0 \geq \lambda_1 \geq \cdots \searrow 0^+$, and a corresponding orthonormal eigenbasis $ \{ \phi_0, \phi_1, \ldots \} $ of $L^2(\rho)$. Moreover,  $ \{ \psi_0, \psi_1, \ldots \} $ with $ \psi_j = K_\rho \phi_j / \lambda_j^{1/2} $ is an orthonormal set in $\mathcal H(M)$. In fact, Mercer's theorem provides an explicit representation $k(\omega,\omega') = \sum_{j=0}^\infty\psi_j(\omega)\psi_j(\omega')$, where direct evaluation of the kernel in the left-hand side (known as ``kernel trick'') avoids the complexity of inner-product computations between feature vectors $ \psi_j$. Here, our perspective is to rely on the orthogonality of the eigenbasis to approximate observables of interest at fixed $L$, and establish convergence of the estimator as $L\to\infty$. A similar approach was adopted for density estimation on non-compact domains, with Mercer-type kernels based on orthogonal polynomials \cite{ZhangEtAl19}.


Now a key operation that the RKHS enables is the \emph{Nystr\"om extension}, which interpolates  $L^2(\rho)$ elements of appropriate regularity to RKHS functions. The Nystr\"om operator $\mathcal N_\rho : D(\mathcal N_\rho) \to \mathcal H(M)$ is defined on the domain $ D(\mathcal N_\rho) = \{ \sum_{j}c_j \phi_j : \sum_j \lvert c_j \rvert^2 / \lambda_j < \infty \}$  by linear extension of $\mathcal N_\rho \phi_j = \psi_j / \lambda_j^{1/2}$. Note that $ \mathcal N_\rho \phi_j = K_\rho \phi_j / \lambda_j = \varphi_j$, so $\mathcal N_\rho$ maps $ \phi_j $ to its continuous representative, and  $K^*_\rho \mathcal N_\rho f = f$, meaning that $ \mathcal N_\rho f =  f$, $\rho$-a.e.  While $D(\mathcal N_\rho)$ may be a strict $L^2(\rho) $ subspace, for any $L$ with $\lambda_{L-1} > 0$ we define a spectrally truncated  operator $ \mathcal N_{L,\rho} : L^2(\rho) \to \mathcal H(M) $, $ \mathcal N_{L,\rho} \sum_j c_j \phi_j = \sum_{j=0}^{L-1} c_j \psi_j / \lambda_j^{1/2}$. Then,  as $ L  $ increases, $K^*_\rho \mathcal N_{L,\rho} f $ converges to $ \Pi_{G_\rho} f $ in $L^2(\rho)$.  To make empirical forecasts, we set $\rho = \mu_N$, compute the expansion coefficients $c_{j,N}(t)$ of  $ U^t Y$ in the $ \{ \phi_{j,N} \} $ basis of $L^2(\mu_N)$, and  construct  $ Y_{t,L,N} = \mathcal N_{L,N} U^t Y \in \mathcal{H}(M) $. Because $ \psi_{j,N} $ are pullbacks of known functions $ u_{j,N} \in C(\mathcal X) $, we have $ Y_{t,L,N} = Z_{t,L,N} \circ X$, where $ Z_{t,L,N} = \sum_{j=0}^{L-1} c_j(t) u_{j,N} / \lambda_{j,N}^{1/2} $ can be evaluated at any $x \in \mathcal X$. 



\begin{figure}
    \centering
    \includegraphics[width=.95\linewidth]{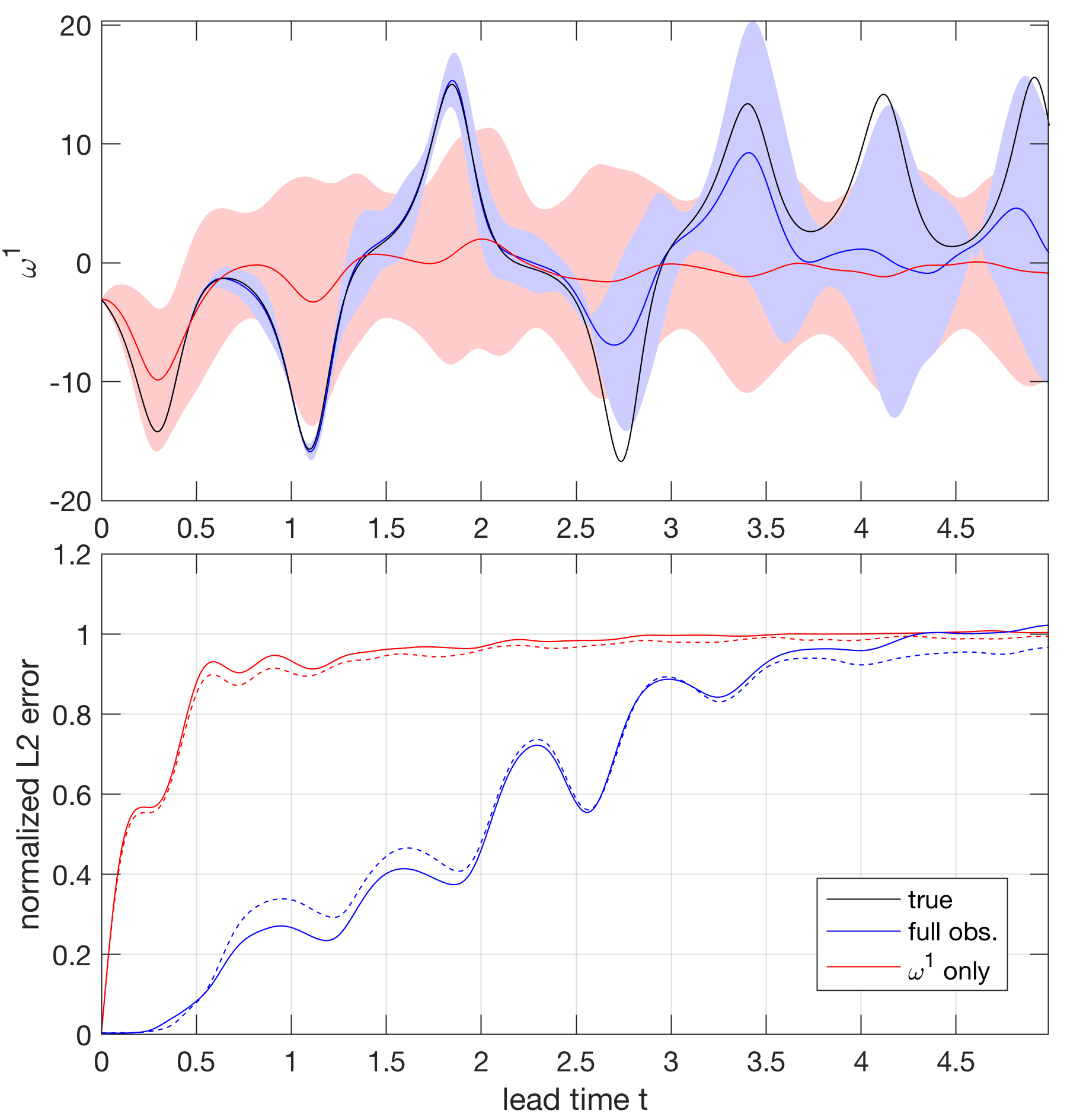}
    \caption{\label{figL63Analog}KAF applied to the L63 state vector component $Y(\omega) = \omega_1$ with full (blue) and partial (red) observations. In the fully observed case, the covariate $X$ is the identity map on  $\Omega = \mathbb R^3$. In the partially observed case, $X(\omega) = \omega_1$ is the projection to the first coordinate. Top: Forecasts $Z_{t,L,N}(x) $ initialized from  fixed $x=X(\omega)$, compared with the true evolution $U^tY(\omega)$ (black). Shaded regions show error bounds based on KAF estimates of the conditional standard deviation, $\sigma_{t}(x)$. Bottom: RMS forecast errors (solid lines) and  $\sigma_{t}$ (dashed lines). The agreement between actual and estimated errors indicates that $\sigma_{t}$ provides useful uncertainty quantification.}
\end{figure}

The function $Y_{t,L,N} $ is our estimator of the conditional expectation $ \mathbb E( U^tY \mid X)$. By spectral convergence of kernel integral operators, as $N\to \infty $, $ Y_{t,L,N}$ converges to $Y_{t,L} := \mathcal N_L U^t Y$ in $C(M)$ norm, where $ \mathcal N_{L} \equiv \mathcal N_{L,\mu}$. Then, as $L\to\infty$, $ K^* Y_{t,L}$ converges in $L^2(\mu)$ norm to $ \Pi_G U^t Y$. Because $\kappa$ is strictly positive-definite, $G$ has dense range in $L^2_X(\mu)$, and thus $\Pi_G U^t Y = \Pi_X U^tY = \mathbb E(U^tY \mid X) $. We therefore conclude that $Y_{t,L,N}$ converges to the conditional expectation as $L \to \infty $ after $N \to\infty$.  Forecast results from the L63 system are shown in Figure~\ref{figL63Analog}. 

\section*{Coherent pattern extraction}

We now turn to the task of coherent pattern extraction in Problem~\ref{probSpectral}. This is a fundamentally unsupervised learning problem, as we seek to discover observables of a dynamical system that exhibit a natural time evolution (by some suitable criterion), rather than approximate a given observable as in the context of forecasting. We have mentioned POD as a technique for identifying coherent observables through eigenfunctions of covariance operators. Kernel PCA \cites{ScholkopfEtAl98} is a generalization of this approach utilizing integral operators with potentially nonlinear kernels. For data lying on Riemannian manifolds, it is popular to employ kernels approximating geometrical operators, such as heat operators and their associated Laplacians. Examples include Laplacian eigenmaps \cite{BelkinNiyogi03}, diffusion maps \cite{CoifmanLafon06}, and variable-bandwidth kernels \cite{BerryHarlim16}. Meanwhile, coherent pattern extraction techniques based on evolution operators have also gained popularity in recent years. These methods include spectral analysis of transfer operators for detection of invariant sets \cites{DellnitzJunge99,DellnitzEtAl00}, harmonic averaging \cite{Mezic05} and dynamic mode decomposition (DMD) \cites{RowleyEtAl09,Schmid10,WilliamsEtAl15,KutzEtAl17,KlusEtAl18} techniques for approximating Koopman eigenfunctions, and Darboux kernels for approximating spectral projectors \cite{KordaEtAl20}. While natural from a theoretical standpoint, evolution operators tend to have more complicated spectral properties than kernel integral operators, including non-isolated eigenvalues and continuous spectrum.  The following examples illustrate distinct behaviors associated with the point ($H_a$) and continuous ($H_c$) spectrum subspaces of $L^2(\mu)$.

\begin{exmpl}[Torus rotation] \label{exmplT2} A quasiperiodic rotation on the 2-torus, $\Omega = \mathbb T^2$, is governed by the system of ODEs $\dot \omega = \vec V(\omega)$, where $\omega = ( \omega_1, \omega_2 ) \in [ 0, 2 \pi )^2$, $ \vec V = (\nu_1,\nu_2)$, and $\nu_1,\nu_2 \in \mathbb R$ are rationally independent frequency parameters. The resulting flow,  $ \Phi^t(\omega) = ( \omega_1 + \nu_1 t, \omega_2 + \nu_2 t ) \mod 2 \pi$, has a unique Borel ergodic invariant probability measure $\mu$ given by a normalized Lebesgue measure. Moreover, there exists an orthonormal basis of $L^2(\mu)$ consisting of Koopman eigenfunctions $z_{jk}(\omega) = e^{i(j\omega_1 + k \omega_2)}$, $j,k \in \mathbb Z$, with eigenfrequencies $\alpha_{jk} = j \nu_1 + k \nu_2$. Thus,  $H_a = L^2(\mu)$, and $H_c$ is the zero subspace. Such a system is said to have a \emph{pure point spectrum}.     
\end{exmpl}

\begin{exmpl}[Lorenz 63 system] \label{exmplL63} The L63 system on $\Omega = \mathbb R^3 $ is governed by a system of smooth ODEs with two quadratic nonlinearities. This system is known to exhibit a physical ergodic invariant probability measure $\mu$ supported on a compact set (the L63 attractor), with mixing dynamics. This means that $H_a$ is the one-dimensional subspace of $L^2(\mu)$ consisting of constant functions, and $H_c$ consists of all $L^2(\mu)$ functions orthogonal to the constants (i.e., with zero expectation value with respect to $\mu$). 
\end{exmpl}

\subsection*{Delay-coordinate approaches}

For the point spectrum subspace $H_a$, a natural class of coherent observables is provided by the Koopman eigenfunctions. Every Koopman eigenfunction $z_j \in H_a $ evolves as a harmonic oscillator at the corresponding eigenfrequency, $ U^t z_{j} = e^{i\alpha_{j}t} z_{j}$, and the associated autocorrelation function, $C_{z_{j}}(t) = e^{i\alpha_{j}t}$, also has a harmonic evolution. Short of temporal invariance (which only occurs for constant eigenfunctions under measure-preserving ergodic dynamics), it is natural to think of a harmonic evolution as being ``maximally'' coherent. In particular, if $z_{j}$ is continuous, then for any $\omega \in A$, the real and imaginary parts of the time series $ t \mapsto U^t z_{j}(\omega) $ are pure sinusoids, even if the flow $\Phi^t$ is aperiodic. Further attractive properties of Koopman eigenfunctions include the facts that they are intrinsic to the dynamical system generating the data, and they are closed under pointwise multiplication, $z_jz_k = z_{j+k}$, allowing one to generate every eigenfunction from a potentially finite generating set. 


Yet, consistently approximating Koopman eigenfunctions from data is a non-trivial task, even for simple systems. For instance, the torus rotation in Example~\ref{exmplT2} has a dense set of eigenfrequencies by rational independence of the basic frequencies $\nu_1$ and $\nu_2$. Thus, any open interval in $\mathbb R$ contains infinitely many eigenfrequencies $\alpha_{jk}$, necessitating some form of regularization. Arguably, the term  ``pure point spectrum'' is somewhat of a misnomer for such systems since a non-empty continuous spectrum is present. Indeed, since the spectrum of an operator on a Banach space includes the closure of the set of eigenvalues,  $i \mathbb R \setminus \{ i \alpha_{jk} \} $ lies in the continuous spectrum.

As a way of addressing these challenges,  observe that if $G$ is a self-adjoint, compact operator commuting with the Koopman group (i.e., $U^t G = G U^t $), then any eigenspace $W_\lambda$ of $G$ corresponding to a nonzero eigenvalue $\lambda$ is invariant under $U^t$, and thus under the generator $V$. Moreover, by compactness of $G$, $W_\lambda$ has finite dimension.  Thus, for any orthonormal basis $ \{ \phi_0, \ldots, \phi_{l-1} \}$ of $W_\lambda $, the generator $V $ on $W_\lambda$ is represented by a skew-symmetric, and thus unitarily diagonalizable, $l \times l $ matrix $ \bm V = [ \langle \phi_i, V \phi_j \rangle_{L^2(\mu)} ] $. The eigenvectors $\vec u = ( u_0, \ldots, u_{l-1} )^\top \in \mathbb C^l$ of $ \bm V $ then contain expansion coefficients of Koopman eigenfunctions $ z = \sum_{j=0}^{l-1} u_j \phi_j$  in $W_\lambda$, and the eigenvalues corresponding to $ \vec u$ are eigenvalues of $ V$. 

On the basis of the above, since any integral operator $G  $ on $L^2(\mu)$ associated with a symmetric kernel $k \in L^2(\mu\times \mu)$ is Hilbert-Schmidt (and thus compact), and we have a wide variety of data-driven tools for approximating integral operators, we can reduce the problem of consistently approximating the point spectrum of the Koopman group on $L^2(\mu)$ to the problem of constructing a commuting integral operator. As we now argue, the success of a number of techniques, including singular spectrum analysis (SSA) \cites{BroomheadKing86,VautardGhil89}, diffusion-mapped delay coordinates (DMDC) \cite{BerryEtAl13}, nonlinear Laplacian spectral analysis (NLSA) \cite{GiannakisMajda12a}, and Hankel DMD \cite{BruntonEtAl17}, in identifying coherent patterns can at least be partly attributed to the fact that they employ integral operators that approximately commute with the Koopman operator. 

A common characteristic of these methods is that they employ, in some form, \emph{delay-coordinate maps} \cite{SauerEtAl91}. With our notation for the covariate function $ X : \Omega \to \mathcal X$ and sampling interval $\Delta t$, the $Q$-step delay-coordinate map is defined as $X_Q : \Omega \to \mathcal X^Q$ with $X_Q(\omega ) = \left( X(\omega_0), X(\omega_{-1}), \ldots,X(\omega_{-Q+1}) \right) $ and $\omega_q = \Phi^{q\,\Delta t}(\omega)$. That is, $X_Q$ can be thought of as a lift of  $X$, which produces ``snapshots'', to a map taking values in the space $\mathcal X^Q$  containing ``videos''. Intuitively, by virtue of its higher-dimensional codomain and dependence on the dynamical flow, a delay-coordinate map such as $X_Q$ should provide additional information about the underlying dynamics on $\Omega$ over the raw covariate map $X$. This intuition has been made precise in a number of ``embedology'' theorems \cite{SauerEtAl91}, which state that under mild assumptions,  for any compact subset $ S \subseteq \Omega$ (including, for our purposes, the invariant set $A$), the delay-coordinate map $X_Q$ is injective on $S$ for sufficiently large $Q$. As a result, delay-coordinate maps provide a powerful tool for state space reconstruction, as well as for constructing informative predictor functions in the context of forecasting. 


Aside from considerations associated with topological reconstruction, however, observe that a metric $d : \mathcal X \times \mathcal X \to \mathbb R $ on covariate space pulls back to a distance-like function $\tilde d_Q : \Omega \times \Omega \to \mathbb R$ such that
\begin{align}
    \label{eqDelayD} \tilde d^2_Q(\omega,\omega') = \frac{1}{Q} \sum_{q=0}^{Q-1}  d^2( X(\omega_{-q} ),X(\omega'_{-q})). 
\end{align}
In particular, $\tilde d_Q^2$ has the structure of an ergodic average of a continuous function under the product dynamical flow $\Phi^t \times \Phi^t$ on $\Omega \times \Omega$. By the von Neumann ergodic theorem, as $Q \to \infty$,  $\tilde d_Q$ converges in $L^2(\mu \times \mu)$ norm to a bounded function $\tilde d_\infty $, which is invariant under the Koopman operator $U^t \otimes U^t$ of the product dynamical system. Note that $\tilde d_\infty $ need not be $\mu \times \mu$-a.e.\ constant, as $\Phi^t \times \Phi^t$ need not be ergodic, and aside from special cases it will not be continuous on $A\times A$. Nevertheless, based on the $L^2(\mu\times\mu)$ convergence of $\tilde d_Q$ to $\tilde d_\infty$, it can be shown \cite{DasGiannakis19} that for any continuous function $h : \mathbb R \to \mathbb R$, the integral operator $G_\infty $ on $  L^2(\mu) $ associated with the kernel $k_\infty = h \circ d_\infty $ commutes with  $U^t $ for any $ t \in \mathbb R$. Moreover, as $Q \to \infty$, the operators $G_Q $ associated with $k_Q = h \circ d_Q$ converge to $G_\infty$ in $L^2(\mu)$ operator norm, and thus in spectrum.    

Many of the operators employed in SSA, DMDC, NLSA, and Hankel DMD can be modeled after $G_Q$ described above. In particular, because $G_Q$ is induced by a continuous kernel, its spectrum can be consistently approximated by data-driven operators $G_{Q,N} $ on $L^2(\mu_N)$, as described in the context of forecasting. The eigenfunctions of these operators at nonzero eigenvalues  approximate eigenfunctions of $G_Q$, which approximate in turn eigenfunctions of $G_\infty$ lying in finite unions of Koopman eigenspaces. Thus, for sufficiently large  $N$ and  $Q$, the eigenfunctions of $G_{Q,N}$ at nonzero eigenvalues capture distinct timescales associated with the point spectrum of the dynamical system, providing physically interpretable features. These kernel eigenfunctions can also be employed in Galerkin schemes to approximate individual Koopman eigenfunctions. 


Besides the spectral considerations described above, in \cite{BerryEtAl13} a geometrical characterization of the eigenspaces of $G_{Q}$ was given based on Lyapunov metrics of dynamical systems. In particular, it follows by Oseledets' multiplicative ergodic theorem that  for $\mu$-a.e.\ $\omega \in M$ there exists a decomposition $T_\omega M = F_{1,\omega}  \oplus \ldots \oplus F_{r,\omega} $, where $ T_\omega M $ is the tangent space at $\omega \in M$, and $F_{j,\omega}$ are subspaces satisfying the equivariance condition $ D\Phi^t F_{j,\omega} = F_{j,\Phi^t(\omega)} $. Moreover, there exist  $ \Lambda_1 >  \cdots > \Lambda_r $, such that for  every $v \in F_{j,\omega} $, $\Lambda_j = \lim_{t\to\infty} \int_0^t \log \lVert D\Phi^s \, v \rVert \, ds / t$, where $ \lVert \cdot \rVert $ is the norm on $T_\omega M$ induced by a Riemannian metric. The numbers $\Lambda_j$ are called \emph{Lyapunov exponents}, and are metric-independent. Note that the dynamical vector field $ \vec V(\omega) $ lies in a subspace $F_{j_0,\omega}$ with a corresponding zero Lyapunov exponent. 

If $F_{j_0,\omega} $ is one-dimensional, and the norms $ \lVert D\Phi^t \, v \rVert$ obey appropriate uniform growth/decay bounds with respect to $ \omega\in M$, the dynamical flow is said to be \emph{uniformly hyperbolic}. If, in addition, the support $A$ of $\mu$ is a differentiable manifold, then  there exists a class of Riemannian metrics, called \emph{Lyapunov metrics}, for which the $F_{j,\omega} $ are mutually orthogonal at every $\omega \in A$. In \cite{BerryEtAl13}, it was shown that using a modification of the delay-distance in~\eqref{eqDelayD} with exponentially decaying weights, as $Q\to\infty $, the top eigenfunctions $ \phi^{(Q)}_{j} $ of $G_Q$ vary predominantly along the subspace $F_{r,\omega}$ associated with the most stable Lyapunov exponent. That is, for every $\omega \in \Omega $ and tangent vector $ v \in T_\omega M $ orthogonal to $F_{r,\omega}$ with respect to a Lyapunov metric, $ \lim_{Q\to\infty} v \cdot \nabla \phi^{(Q)}_j = 0 $.   

\subsection*{RKHS approaches}

While delay-coordinate maps are effective for approximating the point spectrum and associated Koopman eigenfunctions, they do not address the problem of identifying coherent observables in the continuous spectrum subspace $H_c$. Indeed, one can verify that in mixed-spectrum systems the infinite-delay operator $G_\infty $, which provides access to the eigenspaces of the Koopman operator, has a non-trivial nullspace that includes $H_c$ as a subspace  More broadly,  there is no obvious way of identifying coherent observables in $H_c$ as eigenfunctions of an intrinsic evolution operator. As a remedy of this problem, we relax the problem of seeking Koopman eigenfunctions, and consider instead \emph{approximate eigenfunctions}. An observable $ z \in L^2(\mu)$ is said to be an $\epsilon$-approximate eigenfunction of $U^t$ if there exists $\lambda_t \in \mathbb C$ such that  
\begin{equation}
    \label{eqApproxEig}
    \lVert U^t z - \lambda_t z \rVert_{L^2(\mu)} < \epsilon \lVert z \rVert_{L^2(\mu)}. 
\end{equation}
The number $ \lambda_t$ is then said to lie in the $\epsilon$-approximate spectrum of $U^t$. A Koopman eigenfunction is an $\epsilon$-approximate eigenfunction for every $\epsilon > 0 $, so we think of~\eqref{eqApproxEig} as a relaxation of the eigenvalue equation, $U^t z- \lambda_t z = 0$. This suggests that a natural notion of coherence of observables in $L^2(\mu)$, appropriate to both the point and continuous spectrum, is that \eqref{eqApproxEig} holds for $\epsilon \ll 1$ and all $t$ in a ``large'' interval.

We now outline an RKHS-based approach \cite{DasEtAl18}, which identifies observables satisfying this condition through eigenfunctions of a regularized operator $\tilde V_\tau $ on $L^2(\mu)$ approximating $V$ with the properties of (i) being skew-adjoint and compact; and (ii) having eigenfunctions in the domain of the Nystr\"om operator, which maps them to differentiable functions in an RKHS. Here, $\tau $ is a positive regularization parameter such that, as $ \tau \to 0^+$, $ \tilde V_\tau$ converges to $V$ in a suitable spectral sense. We will assume that the forward-invariant, compact manifold $M$ has  $C^1$ regularity, but will not require that the support $A$ of the invariant measure be differentiable.      

With these assumptions, let $k : \Omega \times \Omega \to \mathbb R$ be a symmetric, positive-definite kernel, whose restriction on $M \times M$ is continuously differentiable. Then,  the corresponding RKHS $\mathcal H(M)$ embeds continuously in the Banach space $C^1(M)$ of continuously differentiable functions on $M$, equipped with the standard norm. Moreover, because  $V$ is an extension of the directional derivative $ \vec V \cdot \nabla $ associated with the dynamical vector field, every function in $\mathcal H(M)$ lies, upon inclusion, in $D(V)$. The key point here is that regularity of the kernel induces RKHSs of observables which are guaranteed to lie in the domain of the generator. In particular, the range of the integral operator $G = K^* K $ on $L^2(\mu)$ associated with $k$ lies in $D(V)$, so that $A = VG$ is well-defined. This operator is, in fact, Hilbert-Schmidt, with Hilbert-Schmidt norm bounded by the $C^1(M\times M)$ norm of the kernel $k$.  What is perhaps less obvious is that $G^{1/2} V G^{1/2}$ (which ``distributes'' the smoothing by $G$ to the left and right of $V$), defined on the dense subspace $ \{ f \in L^2(\mu) : G^{1/2} f \in D(V) \} $ is also bounded, and thus has a unique closed extension $ \tilde V: L^2(\mu) \to L^2(\mu)$, which turns out to be Hilbert-Schmidt. Unlike $A$, $ \tilde V $ is skew-adjoint, and thus preserves an important structural property of the generator. By skew-adjointness and compactness of $\tilde V$, there exists an orthonormal basis $ \{ \tilde z_j : j \in \mathbb Z \}  $ of $L^2(\mu)$ consisting of its eigenfunctions $\tilde z_j$, with purely imaginary eigenvalues $i \tilde \alpha_j$.  Moreover, (i) all $ \tilde z_j$ corresponding to nonzero $ \tilde \alpha_j$ lie in the domain of the Nystr\"om operator, and therefore have $C^1$ representatives in $\mathcal H(M)$; and (ii) if $k$ is $L^2(\mu)$-universal, Markov, and ergodic, $\tilde V $ has a simple eigenvalue at zero, in agreement with the analogous property of $V$.

Based on the above, we seek to construct a one-parameter family of such kernels $k_\tau$, with associated RKHSs $\mathcal H_\tau(M)$, such that as $ \tau \to 0^+$, the regularized generators $\tilde V_\tau$ converge to $V$ in a sense suitable for spectral convergence. Here, the relevant notion of convergence is \emph{strong resolvent convergence}; that is, for every element $\lambda$ of the resolvent set of $V$ and every $ f \in L^2(\mu)$, $ ( \tilde V_\tau - \lambda )^{-1} f$ must converge to $(  V - \lambda )^{-1} f $. In that case, for every element $i\alpha$ of the spectrum of $V$ (both point and continuous), there exists a sequence of eigenvalues $i \tilde \alpha_{j_\tau,\tau}$ of $\tilde V_\tau$ converging to $i \alpha$ as $ \tau \to 0^+0$. Moreover, for any $\epsilon > 0 $ and $ T > 0 $, there exists $\tau_* > 0 $ such that for all $ 0 < \tau < \tau_*$ and $ \lvert t \rvert < T$, $ e^{i\alpha_{j_\tau,\tau}t}$ lies in the $\epsilon$-approximate spectrum of $U^t$ and  $\tilde z_{j_\tau,\tau}$ is an $\epsilon$-approximate eigenfunction. 

In \cite{DasEtAl18}, a constructive procedure was proposed for obtaining the kernel family $k_\tau $ through a Markov semigroup on $L^2(\mu)$.  This method has a data-driven implementation, with analogous spectral convergence results for the associated integral operators $G_{\tau,N} $ on $ L^2(\mu_N) $ to those described in the setting of forecasting. Given these operators, we approximate $\tilde V_\tau$ by $\tilde V_{\tau,N} = G_{\tau,N}^{1/2} V_N G_{\tau,N}^{1/2}$, where $V_N $ is a skew-adjoint, \emph{finite-difference approximation} of the generator. For example, $ V_N = ( U^1_N - U^{1*}_N ) / ( 2 \, \Delta t) $ is a second-order finite-difference approximation based on the 1-step shift operator $U^1_N$. See Figure~\ref{figL63} for a graphical representation of a generator matrix for L63.  As with our data-driven approximations of $U^t$, we work with a rank-$L $ operator  $ \hat V_\tau := \Pi_{\tau,N,L}V_{\tau,N} \Pi_{\tau,N,L} $, where $ \Pi_{\tau,N,L} $ is the orthogonal projection to the subspace spanned by the first $L$ eigenfunctions of $G_{\tau,N}$. This family of operators convergences spectrally to $V_\tau$ in a limit of $N \to 0 $, followed by $\Delta t \to 0$ and $L\to \infty$, where we note that $C^1 $ regularity of $k_{\tau} $ is important for the finite-difference approximations to converge.

\begin{figure}
    \centering
    \includegraphics[width=.75\linewidth]{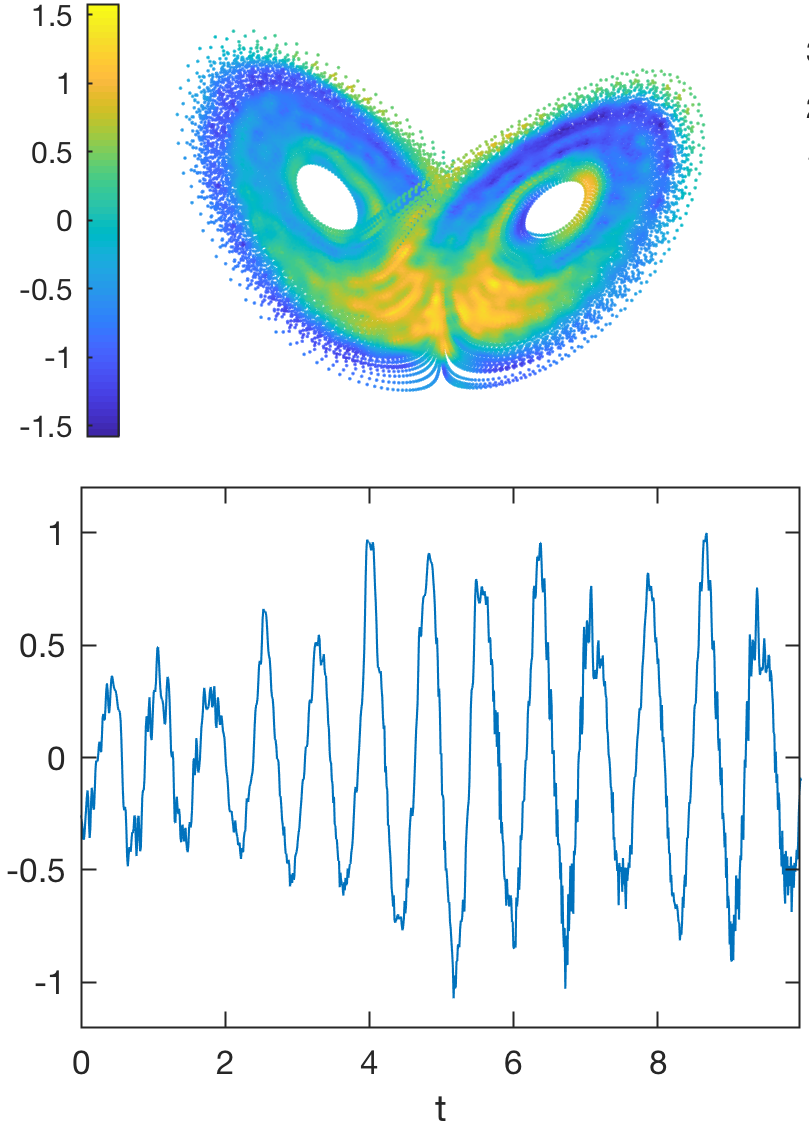}
    \caption{\label{figL63Spec}A representative eigenfunction $\hat z_{j,\tau}$ of the compactified  generator $ \hat V_\tau$ for the L63 system, with low corresponding Dirichlet energy. Top: Scatterplot of $ \Real \hat z_{j,\tau} $ on the L63 attractor. Bottom: Time series of $\Real \hat z_{j,\tau}$ sampled along a dynamical trajectory.} 
\end{figure}

At any given $\tau$, an a posteriori criterion for identifying candidate eigenfunctions $ \hat z_{j,\tau}$ satisfying~\eqref{eqApproxEig} for small $ \epsilon $ is to compute a \emph{Dirichlet energy functional}, $\mathcal D( \hat z_{j,\tau})= \lVert \mathcal N_{\tau,N} \hat z_{j,\tau} \rVert_{\mathcal H_\tau(M)}^2 / \lVert \hat z_{j,\tau} \rVert^2_{L^2(\mu_N)}$. Intuitively, $\mathcal D$ assigns a measure of roughness to every nonzero element in the domain of the Nystr\"om operator (analogously to the Dirichlet energy in Sobolev spaces on differentiable manifolds), and the smaller $\mathcal D(\hat z_{j,\tau})$ is, the more coherent $ \hat z_{j,\tau } $ is expected to be. Indeed, as shown in Figure~\ref{figL63Spec}, the $\hat z_{j,\tau}$ corresponding to low Dirichlet energy identify observables of the L63 system with a coherent dynamical evolution, even though this system is mixing and has no nonconstant Koopman eigenfunctions. Sampled along dynamical trajectories, the approximate Koopman eigenfunctions resemble amplitude-modulated wavepackets, exhibiting a low-frequency modulating envelope while maintaining phase coherence and a precise carrier frequency. This behavior can be thought of as a ``relaxation'' of Koopman eigenfunctions, which generate pure sinusoids with no amplitude modulation.

\section*{Conclusions and outlook}

We have presented mathematical techniques at the interface of dynamical systems theory and data science for statistical analysis and modeling of dynamical systems. One of our primary goals has been to highlight a fruitful interplay of ideas from ergodic theory, functional analysis, and differential geometry, which, coupled with learning theory, provide an effective route for data-driven prediction and pattern extraction, well-adapted to handle nonlinear dynamics and complex geometries. 

There are several open questions and future research directions stemming from these topics. First, it should be possible to combine pointwise estimators derived from methods such as diffusion forecasting and KAF with the Mori-Zwanzig formalism so as to incorporate memory effects. Another potential direction for future development is to incorporate wavelet frames, particularly when the measurements or probability densities are highly localized. Moreover, when the attractor $A$ is not a manifold, appropriate notions of regularity need to be identified so as to fully characterize the behavior of  kernel algorithms such as diffusion maps. While we suspect that kernel-based constructions will still be the fundamental tool, the choice of kernel may need to be adapted to the regularity of the attractor to obtain optimal performance. Finally, a number of applications (e.g., analysis of perturbations) concern the action of dynamics on more general vector bundles besides functions, potentially with a non-commutative algebraic structure, calling for the development of suitable data-driven techniques for such spaces. 

\paragraph*{Acknowledgments} Research of the authors described in this review was supported by DARPA grant HR0011-16-C-0116, NSF grants 1842538, DMS-1317919, DMS-1521775, DMS-1619661, DMS-172317, DMS-1854383, and ONR grants N00014-12-1-0912, N00014-14-1-0150, N00014-13-1-0797, N00014-16-1-2649, N00014-16-1-2888.



\end{document}